\input amstex
\magnification=\magstep1 
\baselineskip=13pt
\documentstyle{amsppt}
\vsize=8.7truein \CenteredTagsOnSplits \NoRunningHeads
\def\tr{\operatorname{tr}}
\def\per{\operatorname{per}}
\def\diag{\operatorname{diag}}
 
 \topmatter

\title  Concentration of the mixed discriminant of well-conditioned matrices  \endtitle
\author Alexander Barvinok  \endauthor
\address Department of Mathematics, University of Michigan, Ann Arbor,
MI 48109-1043, USA \endaddress
\email barvinok$\@$umich.edu \endemail
\date June 2015 \enddate
\thanks  This research was partially supported by NSF Grant DMS 1361541.
\endthanks 
\abstract We call an $n$-tuple $Q_1, \ldots, Q_n$ of positive definite $n \times n$ matrices $\alpha$-conditioned
for some $\alpha \geq 1$ if the ratio of the largest among the eigenvalues of $Q_1, \ldots, Q_n$ to the smallest among the eigenvalues of $Q_1, \ldots, Q_n$ does not exceed $\alpha$. An $n$-tuple is called doubly stochastic if the sum of $Q_i$ is the identity matrix and the trace of each $Q_i$ is 1. We prove that for any fixed $\alpha \geq 1$ the mixed discriminant of
an $\alpha$-conditioned doubly stochastic $n$-tuple is $n^{O(1)} e^{-n}$. As a corollary, for any $\alpha \geq 1$ fixed in advance, we obtain a polynomial time algorithm approximating the mixed discriminant of an $\alpha$-conditioned $n$-tuple within a polynomial in $n$ factor.
\endabstract
\endtopmatter

\document

\head 1. Introduction and main results \endhead

\subhead (1.1) Mixed discriminants \endsubhead Let $Q_1, \ldots, Q_n$ be $n \times n$ real symmetric matrices. The function
$\det \left(t_1 Q_1 + \ldots + t_n Q_n \right)$, where $t_1, \ldots, t_n$ are real variables, is a homogeneous polynomial of degree $n$ in $t_1, \ldots, t_n$ and its coefficient
$$D\left(Q_1, \ldots, Q_n \right) = {\partial^n \over \partial t_1 \cdots \partial t_n} \det\left(t_1 Q_1 +\ldots +t_n Q_n \right)
\tag1.1.1$$ 
is called the {\it mixed discriminant} of $Q_1, \ldots, Q_n$ (sometimes, the normalizing factor of $1/n!$ is used). Mixed discriminants were introduced by A.D. Alexandrov  in his work on mixed volumes \cite{Al38}, see also \cite{Le93}. They also have some interesting combinatorial applications, see Chapter V of \cite{BR97}. 

Mixed discriminants generalize permanents. If the matrices $Q_1, \ldots, Q_n$ are diagonal, so that 
$$Q_i=\diag\left(a_{i1}, \ldots, a_{in} \right) \quad \text{for} \quad i=1, \ldots, n,$$
then 
$$D\left(Q_1, \ldots, Q_n\right) =\per A \quad \text{where} \quad A=\left(a_{ij}\right) \tag1.1.2$$
and 
$$\per A =\sum_{\sigma \in S_n} \prod_{i=1}^n a_{i \sigma(i)}$$ is 
the {\it permanent} of an $n \times n$ matrix $A$. Here the $i$-th row of $A$ is the diagonal of $Q_i$ and $S_n$ is the symmetric  group of all $n!$ permutations of the set $\{1, \ldots, n\}$.

\subhead (1.2) Doubly stochastic $n$-tuples \endsubhead If $Q_1, \ldots, Q_n$ are positive semidefinite matrices then $D\left(Q_1, \ldots, Q_n\right) \geq 0$, see \cite{Le93}. We say that the $n$-tuple $\left(Q_1, \ldots, Q_n \right)$ is 
{\it doubly stochastic} if $Q_1, \ldots, Q_n$ are positive semidefinite,
$$Q_1 + \ldots + Q_n =I \quad \text{and} \quad \tr Q_1 = \ldots =\tr Q_n=1,$$
where $I$ is the $n \times n$ identity matrix and $\tr Q$ is the trace of $Q$. We note that if $Q_1, \ldots, Q_n$ are diagonal
then the $n$-tuple $\left(Q_1, \ldots, Q_n \right)$ is doubly stochastic if and only if the matrix $A$ in (1.1.2) is doubly stochastic, that is, non-negative and has row and column sums 1. 

In \cite{Ba89} Bapat conjectured what should be the mixed discriminant version of the van der Waerden inequality for permanents: if $\left(Q_1, \ldots, Q_n\right)$ is a 
doubly stochastic $n$-tuple then 
$$D\left(Q_1, \ldots, Q_n\right) \ \geq \ {n! \over n^n} \tag1.2.1 $$
where equality holds if and only if 
$$Q_1 = \ldots = Q_n={1 \over n} I.$$
The conjecture was proved by Gurvits \cite{Gu06}, see also \cite{Gu08} for a more general result with a simpler proof.

In this paper, we prove that $D\left(Q_1, \ldots, Q_n\right)$ remains close to $n!/n^n \approx e^{-n}$ if the $n$-tuple 
$\left(Q_1, \ldots, Q_n \right)$ is doubly stochastic and well-conditioned.

\subhead (1.3) $\alpha$-conditioned $n$-tuples \endsubhead For a symmetric matrix $Q$, let 
$\lambda_{\min}(Q)$ denote the minimum eigenvalue of $Q$ and let $\lambda_{\max}(Q)$ denote the maximum eigenvalue of $Q$. We say that a positive definite matrix $Q$ is $\alpha$-{\it conditioned} for some $\alpha \geq 1$ if
$$\lambda_{\max}(Q) \ \leq \ \alpha \lambda_{\min}(Q).$$
We say that an $n$-tuple $\left(Q_1, \ldots, Q_n\right)$ is $\alpha$-{\it conditioned} if 
$$\lambda_{\max}\left(Q_i\right) \ \leq \ \alpha \lambda_{\min}\left(Q_j\right) \quad \text{for all} \quad 1 \leq i, j \leq n. \tag1.3.1$$
In particular, each matrix $Q_i$ is $\alpha$-conditioned, as we allow $i=j$ in (1.3.1).

The main result of this paper is the following inequality.
\proclaim{(1.4) Theorem} Let $\left(Q_1, \ldots, Q_n\right)$ be an $\alpha$-conditioned doubly stochastic $n$-tuple of positive definite $n \times n$ matrices. Then 
$$D\left(Q_1, \ldots Q_n \right) \ \leq \ n^{\alpha^4} e^{-(n-1)}.$$
\endproclaim
Combining the bound of Theorem 1.4 with (1.2.1), we conclude that for any $\alpha \geq 1$, fixed in advance, the mixed discriminant of an $\alpha$-conditioned doubly stochastic $n$-tuple is within a polynomial in $n$ factor of $e^{-n}$. If we allow $\alpha$ to vary with $n$ then as long as $\alpha \ll \root4\of{n \over \ln n}$, the logarithmic order of the mixed discriminant is captured by 
$e^{-n}$.

The estimate of Theorem 1.4 is unlikely to be precise. It can be considered as a (weak) mixed discriminant extension of the Bregman - Minc inequality for permanents (we discuss the connection in Section 1.7).

\subhead (1.5) Scaling \endsubhead We say that an $n$-tuple $\left(P_1, \ldots, P_n\right)$ of $n \times n$ positive definite matrices is obtained from an $n$-tuple $\left(Q_1, \ldots, Q_n\right)$ of $n \times n$ positive definite matrices by {\it scaling}
if for some invertible $n \times n$ matrix $T$ and real $\tau_1, \ldots, \tau_n >0$, we have 
$$P_i =\tau_i T^{\ast} Q_i T \quad \text{for} \quad i=1, \ldots, n,\tag1.5.1$$
where $T^{\ast}$ is the transpose of $T$.
It is easy to check that 
$$D\left(P_1, \ldots, P_n\right) = \left(\det T \right)^2 \left( \prod_{i=1}^n \tau_i \right) D\left(Q_1, \ldots, Q_n \right), \tag1.5.2$$
 provided (1.5.1) holds, see \cite{GS02}.

This notion of scaling extends to $n$-tuples of positive definite matrices the notion of scaling for positive matrices introduced 
by Sinkhorn \cite{Si64}. Gurvits and Samorodnitsky proved in \cite{GS02} that any $n$-tuple of $n \times n$ positive definite matrices can be obtained by scaling from a doubly stochastic $n$-tuple, and, moreover, this can be achieved in polynomial time, as it reduces to solving a convex optimization problem (the gist of their algorithm is given by Theorem 2.1 below).
More generally, Gurvits and Samorodnitsky discuss when an $n$-tuple of positive {\it semi}definite matrices can be scaled to a doubly stochastic $n$-tuple. As is discussed in \cite{GS02}, the inequality (1.2.1), together with the scaling algorithm, the identity (1.5.2) and the inequality 
$$D\left(Q_1, \ldots, Q_n\right) \ \leq \ 1 $$
for doubly stochastic $n$-tuples $\left(Q_1, \ldots, Q_n\right)$, allow one to estimate within a factor of $n!/n^n \approx e^{-n}$ the mixed discriminant of any given $n$-tuple of $n\times n$ positive semidefinite matrices in polynomial time.

In this paper, we prove that if a doubly stochastic $n$-tuple $\left(P_1, \ldots, P_n\right)$ is obtained from an $\alpha$-conditioned $n$-tuple of positive definite matrices then the $n$-tuple $\left(P_1, \ldots, P_n \right)$ is $\alpha^4$-conditioned (see Lemma 2.4 below).
We also prove the following strengthening of Theorem 1.4.

\proclaim{(1.6) Theorem} Suppose that $\left(Q_1, \ldots, Q_n\right)$ is an $\alpha$-conditioned $n$-tuple of $n \times n$ positive definite matrices and suppose that $\left(P_1, \ldots, P_n\right)$ is a doubly stochastic $n$-tuple of positive definite matrices obtained from $\left(Q_1, \ldots, Q_n \right)$ by scaling. Then 
$$D\left(P_1, \ldots, P_n \right) \ \leq \ n^{\alpha^4} e^{-(n-1)}.$$
\endproclaim
Together with the scaling algorithm of \cite{GS02} and the inequality (1.2.1), Theorem 1.6 allows us to approximate in polynomial time the mixed discriminant $D\left(Q_1, \ldots, Q_n\right)$ of an $\alpha$-conditioned $n$-tuple $\left(Q_1, \ldots, Q_n \right)$ within a factor of $n^{\alpha^4}$. Note that the value of $D\left(Q_1, \ldots, Q_n\right)$ may vary within a factor of 
$\alpha^n$.

\subhead (1.7) Connections to the Bregman - Minc inequality \endsubhead The following inequality for permanents of 0-1 matrices was conjectured by Minc \cite{Mi63} and proved by Bregman \cite{Br73}, see also \cite{Sc78} for a much simplified proof: if $A$ is an $n \times n$ matrix with 0-1 entries and row sums $r_1, \ldots, r_n$, then
$$\per A \ \leq \ \prod_{i=1}^n \left(r_i!\right)^{1/r_i}. \tag1.7.1$$
The author learned from A. Samorodnitsky \cite{Sa00} the following restatement of (1.7.1), see also \cite{So03}.
Suppose that $B=\left(b_{ij}\right)$ is an $n \times n$ stochastic matrix (that is, a non-negative matrix with row sums 1) such that 
$$ 0 \ \leq \ b_{ij} \ \leq {1 \over r_i} \quad \text{for all} \quad i, j \tag1.7.2$$
and some positive integers $r_1, \ldots, r_n$.
Then 
$$\per B \ \leq \ \prod_{i=1}^n  {\left(r_i!\right)^{1/r_i} \over r_i}. \tag1.7.3$$
Indeed, the function $B \longmapsto \per B$ is linear in each row and hence its maximum value on the polyhedron of stochastic matrices satisfying (1.7.2) is attained at a vertex of the polyhedron, that is, where $b_{ij} \in \{0, 1/r_i\}$ for all $i,j$. 
Multiplying the $i$-th row of $B$ by $r_i$, we obtain a 0-1 matrix $A$ with row sums $r_1, \ldots, r_n$ and hence (1.7.3) follows by (1.7.1).

Suppose now that $B$ is a doubly stochastic matrix whose entries do not exceed $\alpha/n$ for some $\alpha \geq 1$.
Combining (1.7.3) with the van der Waerden lower bound, we obtain that 
$$\per B = e^{-n} n^{O(\alpha)}. \tag1.7.4$$
Ideally, we would like to obtain a similar to (1.7.4) estimate for the mixed discriminants $D\left(Q_1, \ldots, Q_n\right)$ of doubly stochastic $n$-tuples of positive semidefinite matrices satisfying 
$$\lambda_{\max} \left(Q_i\right) \ \leq \ {\alpha \over n} \quad \text{for} \quad i=1, \ldots, n. \tag1.7.5$$
In Theorem 1.4 such an estimate is obtained under a stronger assumption that the $n$-tuple $\left(Q_1, \ldots, Q_n \right)$ in addition to being doubly stochastic is also $\alpha$-conditioned. This of course implies (1.7.5) but it also prohibits $Q_i$ from having small (in particular, 0) eigenvalues. The question whether a similar to Theorem 1.4 bound can be proven under the the weaker assumption of (1.7.5) together with the assumption that $\left(Q_1, \ldots, Q_n\right)$ is doubly stochastic remains open.

In Section 2 we collect various preliminaries and in Section 3 we prove Theorems 1.4 and 1.6.

\head 2. Preliminaries \endhead

First, we restate a result of Gurvits and Samorodnitsky \cite{GS02} that is at the heart of their algorithm to estimate the mixed discriminant. We state it in the particular case of positive definite matrices.

\proclaim{(2.1) Theorem} Let $Q_1, \ldots, Q_n$ be $n \times n$ positive definite matrices, let 
$H \subset {\Bbb R}^n$ be the hyperplane,
$$H=\left\{ \left(x_1, \ldots, x_n\right): \quad \sum_{i=1}^n x_i =0 \right\}$$ and let
$f: H \longrightarrow {\Bbb R}$ be the function
$$f\left(x_1, \ldots, x_n\right) = \ln \det \left(\sum_{i=1}^n e^{x_i} Q_i \right).$$
Then $f$ is strictly convex on $H$ and attains its minimum on $H$ at a unique point $\left(\xi_1, \ldots, \xi_n\right)$.
Let $S$ be an $n \times n$, necessarily invertible, matrix such that 
$$S^{\ast} S = \sum_{i=1}^n e^{\xi_i} Q_i \tag2.1.1$$
(such a matrix exists since the matrix in the right hand side of (2.1.1) is positive definite).
Let 
$$\tau_i = e^{\xi_i} \quad \text{for} \quad i=1, \ldots, n,$$
let $T=S^{-1}$ and let 
$$B_i =\tau_i T^{\ast}  Q_i T \quad \text{for} \quad i=1, \ldots, n.$$
Then $\left(B_1, \ldots, B_n \right)$ is a doubly stochastic $n$-tuple of positive definite matrices.
\endproclaim

We will need the following simple observation regarding matrices $B_1, \ldots, B_n$ constructed in Theorem 2.1.

\proclaim{(2.2) Lemma} Suppose that for the matrices $Q_1, \ldots, Q_n$ in Theorem 2.1, we have 
$$\sum_{i=1}^n \tr Q_i =n.$$
Then, for the matrices $B_1, \ldots, B_n$ constructed in Theorem 2.1, we have 
$$D\left(B_1, \ldots, B_n \right) \ \geq \ D\left(Q_1, \ldots, Q_n \right).$$
\endproclaim 
\demo{Proof} We have 
$$D\left(B_1, \ldots, B_n\right)=\left( \det T \right)^2 \left(\prod_{i=1}^n \tau_i \right) D\left(Q_1, \ldots, Q_n \right). 
\tag2.2.1$$
Now, 
$$\prod_{i=1}^n \tau_i =\exp\left\{\sum_{i=1}^n \xi_i \right\}=1 \tag2.2.2$$
and 
$$\left(\det T \right)^2 = \left(\det \sum_{i=1}^n e^{\xi_i} Q_i \right)^{-1} =\exp\left\{-f\left(\xi_1, \ldots, \xi_n\right) \right\}. 
\tag2.2.3$$
Since $\left(\xi_1, \ldots, \xi_n\right)$ is the minimum point of $f$ on $H$, we have 
$$f\left(\xi_1, \ldots, \xi_n\right) \ \leq \ f(0, \ldots, 0) = \ln \det Q \quad \text{where} \quad Q=\sum_{i=1}^n Q_i. \tag2.2.4$$
We observe that $Q$ is a positive definite matrix with eigenvalues, say, $\lambda_1, \ldots, \lambda_n$ such that 
$$\sum_{i=1}^n \lambda_i = \tr Q =\sum_{i=1}^n \tr Q_i =n \quad \text{and} \quad \lambda_1, \ldots, \lambda_n >0.$$
Applying the arithmetic - geometric mean inequality, we obtain
$$\det Q =\lambda_1 \cdots \lambda_n \ \leq \ \left({\lambda_1 + \ldots + \lambda_n \over n} \right)^n \ \leq \ 1. \tag2.2.5$$
Combining (2.2.1) -- (2.2.5), we complete the proof.
{\hfill \hfill \hfill} \qed
\enddemo

\subhead (2.3) From symmetric matrices to quadratic forms \endsubhead With an $n \times n$ symmetric matrix $Q$ we associate the quadratic form $q: {\Bbb R}^n \longrightarrow {\Bbb R}$ defined by 
$$q(x) =\langle Qx, x \rangle \quad \text{for} \quad x \in {\Bbb R}^n,$$
where $\langle \cdot, \cdot \rangle$ is the standard inner product in ${\Bbb R}^n$. We define the eigenvalues, the trace, and the determinant of $q$ as those of $Q$. Consequently, we define the mixed discriminant $D\left(q_1, \ldots, q_n\right)$ of 
quadratic forms $q_1, \ldots, q_n$. An $n$-tuple of positive semidefinite quadratic forms 
$q_1, \ldots, q_n: {\Bbb R}^n \longrightarrow {\Bbb R}$ is {\it doubly stochastic} if 
$$\sum_{i=1}^n q_i(x)=\|x\|^2 \quad \text{for all} \quad x \in {\Bbb R}^n \quad \text{and} \quad 
\tr q_1 =\ldots = \tr q_n =1.$$
The property of being $\alpha$-conditioned extends to positive definite quadratic forms in a natural way. Namely, 
a positive definite quadratic form if $\alpha$-{\it conditioned}, if 
$$q(x) \ \leq \ \alpha q(y) \quad \text{for any two} \quad x, y \in {\Bbb R}^n \quad \text{such that} \quad \|x\|=\|y\|=1,$$
where $\| \cdot\|$ is the standard Euclidean norm in ${\Bbb R}^n$. Similarly, an $n$-tuple of positive definite 
quadratic forms $q_1, \ldots, q_n: {\Bbb R}^n \longrightarrow {\Bbb R}$ is $\alpha$-{\it conditioned}, if each form 
$q_i$ is $\alpha$-conditioned and if 
$$q_i(x) \ \leq \ \alpha q_j(x) \quad \text{for all} \quad x \in {\Bbb R}^n \quad \text{and all} \quad 1 \leq i, j \leq n.$$
An $n$-tuple of quadratic forms $p_1, \ldots, p_n: {\Bbb R}^n \longrightarrow {\Bbb R}^n$ is obtained from an 
$n$-tuple $q_1, \ldots, q_n: {\Bbb R}^n \longrightarrow {\Bbb R}$ by {\it scaling} if for some invertible linear transformation 
$T: {\Bbb R}^n \longrightarrow {\Bbb R}^n$ and real $\tau_1, \ldots, \tau_n > 0$ we have
$$p_i(x)=\tau_i q_i(Tx) \quad \text{for all} \quad x \in {\Bbb R}^n \quad \text{and all} \quad i=1, \ldots, n.$$
One advantage of working with quadratic forms as opposed to matrices is that it is particularly easy to define the restriction of a
quadratic form onto a subspace. We will use the following construction: suppose that 
$q_1, \ldots, q_n: {\Bbb R}^n \longrightarrow {\Bbb R}$ are positive definite quadratic forms and let 
$L \subset {\Bbb R}^n$ be an $m$-dimensional subspace for some $1 \leq m \leq n$. Then $L$ inherits Euclidean structure from ${\Bbb R}^n$ and we can consider the {\it restrictions} $\widehat{q}_1, \ldots, \widehat{q}_n: L \longrightarrow {\Bbb R}$ of $q_1, \ldots, q_n$ onto $L$. Thus we can define the mixed discriminant $D\left(\widehat{q}_1, \ldots, \widehat{q}_m\right)$.
Note that by choosing an orthonormal basis in $L$, we can associate $m \times m$ symmetric matrices 
$\widehat{Q}_1, \ldots, \widehat{Q}_m$ with $\widehat{q}_1, \ldots, \widehat{q}_m$. A different choice of an orthonormal basis results in the transformation $\widehat{Q}_i \longmapsto U^{\ast} \widehat{Q}_i U$ for some $m \times m$ orthogonal matrix $U$ and $i=1, \ldots, m$, which does not change the mixed discriminant $D\left(\widehat{Q}_1, \ldots, \widehat{Q}_m\right)$.

\proclaim{(2.4) Lemma} 
Let $q_1, \ldots, q_n: {\Bbb R}^n \longrightarrow {\Bbb R}$ be an $\alpha$-conditioned $n$-tuple of positive definite quadratic forms. Let $L \subset {\Bbb R}^n$ be an $m$-dimensional subspace, where $1 \leq m \leq n$, let
$T: L \longrightarrow {\Bbb R}^n$ be a linear transformation such that $\ker T=\{0\}$ and let $\tau_1, \ldots, \tau_m >0$
be reals. Let us define quadratic forms $p_1, \ldots, p_m: L \longrightarrow {\Bbb R}$ by
$$p_i(x)=\tau_i q_i(Tx) \quad \text{for} \quad x \in L \quad \text{and} \quad i=1, \ldots, m.$$
Suppose that 
$$\sum_{i=1}^m p_i(x) =\|x\|^2 \quad \text{for all} \quad x \in L \quad  \text{and} \quad \tr p_i=1 \quad \text{for} \quad i=1, \ldots, m.$$
Then the $m$-tuple of quadratic forms $p_1, \ldots, p_m$ is $\alpha^4$-conditioned.
\endproclaim
\demo{Proof} Since the $n$-tuple $q_1, \ldots, q_n$ is $\alpha$-conditioned, we have 
$$q_i(x) \ \leq \ \alpha q_j(x) \quad \text{for all} \quad x \in {\Bbb R}^n \quad \text{and all} \quad 1 \leq i, j \leq n.$$
We define quadratic forms $r_i: L \longrightarrow {\Bbb R}$, $i=1, \ldots, m$, by
$$r_i(x)=q_i(Tx) \quad \text{for} \quad x \in L \quad \text{and} \quad i=1, \ldots, m.$$
Then
$$ r_i(x) \ \leq \ \alpha r_j(x) \quad \text{for all} \quad 1 \leq i, j \leq m \quad \text{and all} \quad x \in L. \tag2.4.1$$
Therefore,
$$\tr r_i \ \leq \ \alpha \tr r_j \quad \text{for all} \quad 1\leq i, j \leq m.$$
Since $1=\tr p_i =\tau_i \tr r_i$, we conclude that $\tau_i =1/\tr r_i$ and, therefore, 
$$\tau_i \ \leq \ \alpha \tau_j \quad \text{for all} \quad 1\leq i, j \leq m \tag2.4.2$$ 
Since $p_i=\tau_i r_i$, combining (2.4.1) and (2.4.2), we obtain
$$p_i(x) \ \leq \ \alpha^2 p_j(x) \quad \text{for all} \quad x \in L \quad \text{and all} \quad 1 \leq i, j \leq m. \tag2.4.3$$
Seeking a contradiction, suppose that 
$$\split p_j(x) \ > \ \alpha^4 p_j(y) \quad &\text{for some} \quad x, y \in L \quad \text{such that} \quad \|x\|=\|y\|=1 \\ &\text{and some} \quad 1 \leq j \leq m. \endsplit$$
Then, applying (2.4.3) twice, we obtain
$$p_i(y) \ \leq \ \alpha^2 p_j(y) \ < \ \alpha^{-2} p_j(x) \ \leq \ p_i(x) \quad \text{for all} \quad i,$$
so in the end 
$$\split p_i(y) \ < \ p_i(x) \quad &\text{for some} \quad x, y \in L \quad \text{such that} \quad \|x\|=\|y\|=1 \\ &\text{and all} \quad i=1,\ldots, m, \endsplit$$
which is a contradiction since 
$$1 = \sum_{i=1}^m p_i(y) = \sum_{i=1}^m p_i(x).$$
This proves that 
$$\split p_j(x) \ \leq \ \alpha^4 p_j(y) \quad &\text{for all} \quad x, y \in L \quad \text{such that} \quad \|x\|=\|y\|=1 \\ &\text{and all} \quad 1 \leq j \leq m \endsplit$$
and hence concludes the proof.
{\hfill \hfill \hfill} \qed
\enddemo

\proclaim{(2.5) Lemma} Let $q_1, \ldots, q_n: {\Bbb R}^n \longrightarrow {\Bbb R}$ be positive semidefinite quadratic forms and suppose that 
$$q_n(x) = \langle u, x \rangle^2,$$
where $u \in {\Bbb R}^n$ and $\|u\|=1$. Let $H=u^{\bot}$ be the orthogonal complement to $u$. Let 
$\widehat{q}_1, \ldots, \widehat{q}_{n-1}: H \longrightarrow {\Bbb R}$ be the restrictions of $q_1, \ldots, q_{n-1}$ onto $H$. Then
$$D(q_1, \ldots, q_n) = D\left(\widehat{q}_1, \ldots, \widehat{q}_{n-1}\right).$$
\endproclaim
\demo{Proof} Let us choose an orthonormal basis of ${\Bbb R}^n$ for which $u$ is the last basis vector and let $Q_1, \ldots, Q_n$ be the matrices of the forms $q_1, \ldots, q_n$ in that basis. Then the only non-zero entry of $Q_n$ is 1 in the lower right corner. Let $\widehat{Q}_1, \ldots, \widehat{Q}_{n-1}$ be the upper left $(n-1) \times (n-1)$ submatrices of 
$Q_1, \ldots, Q_{n-1}$.
Then 
$$\det \left(t_1 Q_1 + \ldots + t_n Q_n \right) = t_n \det \left(t_1 \widehat{Q}_1 + \ldots + t_{n-1} \widehat{Q}_{n-1} \right)$$
and hence by (1.1.1) we have 
$$D\left(Q_1, \ldots, Q_n\right) = D\left(\widehat{Q}_1, \ldots, \widehat{Q}_{n-1}\right).$$
On the other hand, $\widehat{Q}_1, \ldots, \widehat{Q}_{n-1}$ are the matrices of $\widehat{q}_1, \ldots, \widehat{q}_{n-1}$.
{\hfill \hfill \hfill} \qed
\enddemo

Finally, the last lemma before we embark on the proof of Theorems 1.4 and 1.6.

\proclaim{(2.6) Lemma} Let $q: {\Bbb R}^n \longrightarrow {\Bbb R}$ be an $\alpha$-balanced quadratic form such that 
$\tr q=1$. Let $H \subset {\Bbb R}^n$ be a hyperplane and let $\widehat{q}$ be the restriction of $q$ onto $H$. 
Then 
$$1- {\alpha \over n} \ \leq \  \tr \widehat{q} \ \leq \ 1-{1 \over \alpha n}.$$
\endproclaim 
\demo{Proof} Let 
$$0 \ < \ \lambda_1 \ \leq \ \ldots \ \leq \ \lambda_n$$
be the eigenvalues of $q$. Then
$$\sum_{i=1}^n \lambda_i=1 \quad \text{and} \quad  \lambda_n \ \leq \ \alpha \lambda_1,$$
from which it follows that 
$$\lambda_1 \ \geq \ {1 \over \alpha n} \quad \text{and} \quad \lambda_n \ \leq \ {\alpha \over  n}.$$
As is known, the eigenvalues of $\widehat{q}$ interlace the eigenvalues of $q$, see, for example, Section 1.3 of \cite{Ta12}, so for the eigenvalues $\mu_1, \ldots, \mu_{n-1}$ of $\widehat{q}$ we have
$$\lambda_1\ \leq \ \mu_1 \ \leq \ \lambda_2 \ \leq \ \ldots \ \leq \ \lambda_{n-1} \ \leq \ \mu_{n-1} \ \leq \ \lambda_n.$$
Therefore,
$$1-{1 \over \alpha n} \ \geq \ \sum_{i=2}^n \lambda_i \ \geq\ \tr \widehat{q} = \sum_{i=1}^{n-1} \mu_i \ \geq \ \sum_{i=1}^{n-1} \lambda_i \ \geq \ 1 - {\alpha \over  n}.$$
{\hfill \hfill \hfill} \qed
\enddemo

\head 3. Proof of Theorem 1.4 and Theorem 1.6  \endhead

Clearly, Theorem 1.6 implies Theorem 1.4, so it suffices to prove the former.

\subhead (3.1) Proof of Theorem 1.6 \endsubhead As in Section 2.3, we associate quadratic forms with matrices. We prove the following statement by induction on $m=1, \ldots, n$.
\bigskip
\noindent{\bf Statement:}  Let $q_1, \ldots, q_n: {\Bbb R}^n \longrightarrow {\Bbb R}$ be an $\alpha$-conditioned
$n$-tuple of quadratic forms. Let $L \subset {\Bbb R}^n$ be an $m$-dimensional subspace, $1 \leq m \leq n$, let $T: L \longrightarrow {\Bbb R}^n$ be a linear transformation such that $\ker T=\{0\}$ and let $\tau_1, \ldots, \tau_m>0$ be reals. Let us define quadratic forms 
$p_i: L \longrightarrow {\Bbb R}$, $i=1, \ldots, m$, by
$$p_i(x) =\tau_i q_i(Tx) \quad \text{for} \quad x \in L \quad \text{and} \quad i=1, \ldots, m$$
and suppose that
$$\sum_{i=1}^m p_i(x) =\|x\|^2 \quad \text{for all} \quad x \in L \quad  \text{and} \quad \tr p_i=1 \quad \text{for} \quad i=1, \ldots, m.$$
Then 
$$D(p_1, \ldots, p_m) \ \leq \ \exp\left\{-(m-1) + \alpha^4 \sum_{k=2}^m {1 \over k} \right\}. \tag3.1.1$$
\bigskip
In the case of $m=n$, we get the desired result.

The statement holds if $m=1$ since in that case $D(p_1)=\det p_1 =1$.

 Suppose that $m >1$.
Let $L \subset {\Bbb R}^n$ be an $m$-dimensional subspace and let the linear transformation $T$, numbers $\tau_i$ and the forms $p_i$ for $i=1, \ldots, m$ be as above.
By Lemma 2.4, the $m$-tuple 
$p_1, \ldots, p_m$ is $\alpha^4$-conditioned. We write the spectral decomposition
$$p_m(x) =\sum_{j=1}^m \lambda_j \langle u_j, x \rangle^2,$$
where $u_1, \ldots, u_m \in L$ are the unit eigenvectors of $p_m$ and $\lambda_1, \ldots, \lambda_m >0$ are the corresponding eigenvalues of $p_m$. Since $\tr p_m=1$, we have $\lambda_1 + \ldots +\lambda_m =1$.
Let $L_j =u_j^{\bot}$, $L_j \subset L$, be the orthogonal complement of $u_j$ in $L$. Let 
$$\widehat{p}_{ij}:\ L_j \longrightarrow {\Bbb R} \quad \text{for} \quad i=1, \ldots, m-1 \quad \text{and} \quad j=1, \ldots, m$$
be the restriction of $p_i$ onto $L_j$.

Using Lemma 2.5, we write
$$\aligned D(p_1, \ldots, p_m) =& \sum_{j=1}^m \lambda_j D\left(p_1, \ldots, p_{m-1}, \langle u_j, x \rangle^2\right)\\=
&\sum_{j=1}^m \lambda_j D\left(\widehat{p}_{1j}, \ldots, \widehat{p}_{(m-1)j}\right) \quad \text{where} \\
&\sum_{j=1}^m \lambda_j=1 \quad \text{and} \quad \lambda_j >0 \quad \text{for} \quad 
j=1, \ldots, m. \endaligned \tag3.1.2$$
Let
$$\sigma_j = \tr \widehat{p}_{1j} + \ldots + \tr \widehat{p}_{(m-1)j} \quad \text{for} \quad j=1, \ldots, m.$$
Since
$$\sum_{i=1}^{m-1} \widehat{p}_{ij}(x)=\|x\|^2 - p_{mj}(x) \quad \text{for all} \quad x \in L_j \quad \text{and} \quad 
j=1, \ldots, m$$
and since the form $p_{mj}$ is $\alpha^4$-balanced, by Lemma 2.6, we have 
$$m-2 +{1 \over \alpha^4 m} \ \leq \ \sigma_j \ \leq \ m-2 +{\alpha^4 \over m} \quad \text{for} \quad j=1, \ldots, m. \tag3.1.3$$
Let us define
$$r_{ij} = {m-1 \over \sigma_j} \widehat{p}_{ij} \quad \text{for} \quad i=1, \ldots, m-1 \quad \text{and} \quad j=1, \ldots, m.$$
Then by (3.1.3), 
$$\aligned &D\left(\widehat{p}_{1j}, \ldots, \widehat{p}_{(m-1)j}\right) = \left( {\sigma_j \over m-1}\right)^{m-1} 
D\left(r_{1j}, \ldots, r_{(m-1)j}\right) \\ &\quad \leq \ \left(1 -{1 \over m-1} + {\alpha^4 \over m(m-1)}\right)^{m-1} 
D\left(r_{1j}, \ldots, r_{(m-1)j}\right) \\ &\quad \leq \ \exp\left\{ -1 +{\alpha^4 \over m}\right\} D\left(r_{1j}, \ldots, r_{(m-1)j}\right) \\
&\qquad \qquad \text{for} \quad j=1, \ldots, m.
\endaligned \tag3.1.4$$
In addition, 
$$\tr r_{1j} + \ldots + \tr r_{(m-1)j} = m-1 \quad \text{for} \quad j=1, \ldots, m. \tag3.1.5$$
For each $j=1, \ldots, m$, let 
$w_{1j}, \ldots, w_{(m-1)j}:\ L_j \longrightarrow {\Bbb R}$ be a doubly stochastic $(m-1)$-tuple of quadratic forms obtained from 
$r_{1j}, \ldots, r_{(m-1)j}$ by scaling as described in Theorem 2.1. From (3.1.5) and Lemma 2.2, we have 
$$D\left(r_{1j}, \ldots, r_{(m-1)j}\right) \ \leq \ D\left(w_{1j}, \ldots, w_{(m-1)j}\right) \quad \text{for} \quad j=1, \ldots, m.
\tag3.1.6$$
Finally, for each $j=1, \ldots, m$, we are going to apply the induction hypothesis to the $(m-1)$-tuple of quadratic forms $w_{1j}, \ldots, w_{(m-1)j}: L_j \longrightarrow {\Bbb R}$. Since the $(m-1)$-tuple is doubly stochastic, we have 
$$\aligned &\sum_{i=1}^{m-1} w_{ij}(x) =\|x\|^2 \quad \text{for all} \quad x \in L_j \quad \text{and all} \quad j=1, \ldots, m \\
& \qquad \qquad \qquad \text{and} \\
&\tr w_{ij}=1 \quad \text{for all} \quad i=1, \ldots, m-1 \quad \text{and} \quad  j=1, \ldots, m. \endaligned \tag3.1.7$$
Since the $(m-1)$-tuple $w_{1j}, \ldots, w_{(m-1)j}$ is obtained from the $(m-1)$-tuple \break
$r_{1j}, \ldots, r_{(m-1)j}$ by scaling, there are invertible linear operators $S_j: L_j \longrightarrow L_j$ and real numbers 
$\mu_{ij} >0$ for $i=1, \ldots, m-1$ and $j=1, \ldots, m$ such that 
$$\split w_{ij}(x) = \mu_{ij} r_{ij}(S_j x) \quad &\text{for all} \quad x \in L_j \\ & \text{and all} \quad i=1, \ldots, m-1 \quad 
\text{and} \quad j=1, \ldots, m. \endsplit$$
In other words, 
$$\aligned w_{ij}(x) = &\mu_{ij} r_{ij}\left(S_j x\right) ={\mu_{ij}(m-1) \over \sigma_j} \widehat{p}_{ij}\left(S_jx\right)={\mu_{ij}(m-1) \over \sigma_j} p_i \left(S_jx\right)\\=&{\mu_{ij}(m-1) \tau_i \over \sigma_j} q_i\left(TS_jx\right) \quad \text{for all } \quad x \in L_j \\
&\qquad \qquad \text{and all} \quad i=1, \ldots, m-1 \quad \text{and} \quad j=1, \ldots, m.
 \endaligned \tag3.1.8$$
 Since for each $j=1, \ldots, m$, the linear transformation $TS_j: L_j \longrightarrow {\Bbb R}^n$ of an $(m-1)$-dimensional subspace $L_j \subset {\Bbb R}^n$ has zero kernel, from (3.1.7) and (3.1.8) we can apply the induction hypothesis to conclude that 
 $$\aligned &D\left(w_{1j}, \ldots, w_{(m-1)j}\right) \ \leq \ \exp\left\{-(m-2) +\alpha^4 \sum_{k=2}^{m-1} {1 \over k} \right\} \\
& \qquad  \qquad \text{for} \quad j=1, \ldots, m \endaligned
 \tag3.1.9$$
 Combining the inequalities (3.1.2), (3.1.4), (3.1.6) and (3.1.9), we obtain (3.1.1) and conclude the induction step.
{\hfill \hfill \hfill} \qed

\enddocument
\end